\documentclass[11pt]{amsart}
\usepackage{amssymb,amsmath,amsthm}
\usepackage{epsfig}
\newtheorem{theorem}{Theorem}[section]
\newtheorem{lemma}[theorem]{Lemma}
\newtheorem{remark}[theorem]{Remark}
\newtheorem{definition}[theorem]{Definition}

\newtheorem{corollary}[theorem]{Corollary}

\newtheorem{example}[theorem]{Example}

\newcommand\A{{\mathcal A}}
\newcommand\B{{\mathcal B}}
\newcommand\M{{\mathbb M}}
\title{When the theories meet:\ 
Khovanov homology as Hochschild homology of links}
\author{Jozef H. Przytycki}
\begin{document}
\maketitle
\centerline{Version of October 6, 2005}

\begin{quotation}
ABSTRACT.
\baselineskip=10pt
We show that Khovanov homology and Hochschild homology theories share 
common structure. In fact they overlap: Khovanov homology of a 
$(2,n)$-torus link can be interpreted as a Hochschild homology of 
the algebra underlining the Khovanov homology. In the classical case 
of Khovanov homology we prove the concrete connection. In the general case 
of Khovanov-Rozansky, $sl(n)$, homology and their deformations we 
conjecture the connection. The best framework to explore our ideas is to use 
a comultiplication-free version of Khovanov homology for graphs
 developed by L. Helme-Guizon and Y. Rong and extended here
 to $\mathbb M$-reduced case, and in the case of a polygon to 
noncommutative algebras. 
In this framework we prove that for any unital algebra $\A$ 
the Hochschild homology of $\A$ is isomorphic 
to graph homology over $\A$ of a polygon. We expect that this paper will 
encourage a flow of ideas in both directions between Hochschild/cyclic 
homology and Khovanov homology theories.  
\\
\end{quotation}
\ \\

\section{Hochschild homology and cyclic homology}\label{1}
\markboth{\hfil{\sc When the theories meet }\hfil}
{\hfil{\sc Khovanov homology as Hochschild homology of links}\hfil}

We recall in this section definition of Hochschild homology and cyclic 
homology and we sketch two classical calculations for tensor algebras 
and symmetric tensor algebras.
More calculations are  reviewed in Section 4 in which we use our main 
result, Theorem 3.1, to obtain new results in Khovanov homology, 
in particular solving some conjectures from \cite{H-P-R}.\ 
 We follow \cite{Lo} in our exposition of Hochschild homology.

Let $k$ be a commutative ring and $\A$ a $k$-algebra (not necessarily 
commutative). 
Let $\M$ be a bimodule over $\A$ that is a $k$-module on 
which $\A$ operates linearly on the left and on the right in such a way 
that $(am)a'=a(ma')$ for $a,a'\in \A$ and $m \in \M$. The actions of $\A$ 
and $k$ are always compatible (e.g. 
$m(\lambda a)= (m\lambda)a=\lambda (ma)$). When $\A$ has a unit element 
$1$ we always assume that $1m=m1=m$ for all $m \in \M$. Under this 
unital hypothesis, the bimodule $\M$ is equivalent to a right 
$\A\otimes \A^{op}$-module via $m(a'\otimes a)=ama'$. Here $\A^{op}$ 
denotes the opposite algebra of $\A$ that is $\A$ and $\A^{op}$ are 
the same as sets but the product $a\cdot b$ in $\A^{op}$ is the 
product $ba$ in $\A$. The product map of $\A$ is usually denoted 
$\mu: \A \otimes \A \to A$,\ $\mu (a,b)=ab$.\\  

In this paper we work only with unital algebras. We also assume, 
unless otherwise stated, that $\A$ is a free $k$-module, however 
in most cases, it suffices to assume that $\A$ is $k$-projective, or less 
restrictively, that $\A$ is flat over $k$. \ Throughout the paper the 
tensor product $\A \otimes \B$ denotes the tensor product over $k$, 
that is $\A \otimes_k \B$. 

\begin{definition}[\cite{Hoch,Lo}]\label{1.1}\ 
The Hochschild chain complex $C_*(\A,\M)$ is defined as:
$$ \ldots
\stackrel{b}{\to} \M\otimes \A^{\otimes n} \stackrel{b}{\to} 
\M\otimes \A^{\otimes n-1} \stackrel{b}{\to} \ldots 
\stackrel{b}{\to} \M\otimes \A \stackrel{b}{\to} \M$$
where $C_n(\A,\M)= \M\otimes \A^{\otimes n}$ and the Hochschild 
boundary is the $k$-linear map $b: \M\otimes \A^{\otimes n} \to 
\M\otimes \A^{\otimes n-1}$ given by the formula 
$b= \sum_{i=0}^n (-1)^id_i$, where the face maps $d_i$ are given by\\
$d_0(m,a_1,\ldots, a_n)= (ma_1,a_2,\ldots, a_n),$\\
$d_i(m,a_1,\ldots, a_n)= (m,a_1,\ldots , a_ia_{i+1},\ldots , a_n)$ for 
$1\leq i \leq n-1$, \\
$d_n(m,a_1,\ldots, a_n)= (a_nm,a_1,\ldots ,a_{n-1})$.\\
In the case when $\M=\A$ the Hochschild complex is called the 
{\it cyclic bar complex}.\\
By definition, the $n$th Hochschild homology group of the unital 
$k$-algebra $\A$ with coefficients in the $\A$-bimodule $\M$ is 
 the $n$th homology group of the Hochschild chain complex 
denoted by $H_n(\A,\M)$. In the particular case $\M=\A$ we write 
$C_*(\A)$ instead of $C_*(\A,\A)$ and $HH_*(\A)$ instead of $H_*(\A,\A)$.
\end{definition}

The algebra $\A$ acts on $C_n(\A,\M)$ by 
$a\cdot (m,a_1,...,a_n)=(am,a_1,...,a_n)$. 
If $\A$ is a commutative algebra then the action commutes with 
boundary map $b$, therefore $H_n(\A,\M)$ (in particular, $HH_*(\A)$) is 
an $\A$-module. 

If $\A$ is a graded algebra and $\M$ a coherently graded 
$\A$-bimodule, and the boundary maps are grading preserving, then 
the Hochschild chain complex is a bigraded chain complex with 
($b: C_{i,j}(\A,\mathbb M) \to C_{i-1,j}(\A,\mathbb M))$, and 
$H_{**}(\A,\mathbb M)$ is a bigraded $k$-module. In the case of 
abelian $\A$ and $\A$-symmetric $\mathbb M$ (i.e. $am=ma$), 
$H_{**}(\A,\mathbb M)$ is bigraded $\A$-module. The main examples 
coming from the knot theory are $\A_m=Z[x]/(x^m)$ and $\mathbb M$ the 
ideal in $\A_m$ generated by $x^{m-1}$.

We complete this survey by describing, after \cite{Lo},
 two classical results in 
Hochschild homology -- the computation of Hochschild homology for 
a tensor algebra and for a symmetric tensor algebra.
\begin{theorem}\label{1.2}\ \\
Let $V$ be any $k$-module and let $\A =T(V)= k \oplus V 
\oplus V^{\otimes 2} \oplus ...$ be its tensor algebra. 
We denote by $\tau_n : V^{\otimes n} \to V^{\otimes n}$ the cyclic 
permutation, $\tau_n(v_1,...,v_{n-1},v_n) = (v_n,v_1,...,v_{n-1})$.
Then the Hochschild homology of $\A=T(V)$ is:\\ 
$HH_0(\A)= \oplus_{i\geq 0} V^{\otimes i}/(1-\tau_i)  $,\\
$HH_1(\A)= \oplus_{i\geq 1} (V^{\otimes i})^{\tau_i}$, where 
$(V^{\otimes i})^{\tau_i}$ is the space of invariants, that is 
the kernel of $1-\tau_i$.\\
$HH_n(\A)= 0$ for $n \geq 2$.
\end{theorem}
The main idea of the proof is to show that the Hochschild chain complex 
of $T(V)$ is quasi-isomorphic\footnote{Two chain complexes $C$ and $C'$ are 
called quasi-isomorphic if there is a chain map $f: C \to C'$ which 
induces an isomorphism on homology, $f_*: H_n(C) \to H_n(C')$, for all $n$.
The map $f$ is called a quasi-isomorphism, chain equivalence, or homologism.}
 to the ``small" chain complex:
$$ C^{small}(T(V)):\ \ \ ... \to 0 \to \A\otimes V \stackrel{\hat b}{\to} \A $$ 
where the module $\A$ is in degree $0$ and where the map $\hat b$ is 
given by $\hat b(a\otimes v)= av - va$. Therefore $\hat b$ restricted 
to $V^{\otimes n-1}\otimes V$ is precisely 
$(1-\tau_n): V^{\otimes n} \to V^{\otimes n}$ and Theorem 1.2 follows.
 See Proposition 3.1.2 and Theorem 3.1.4 of \cite{Lo}.

\begin{theorem}\label{1.3} Symmetric and polynomial algebras. \\
 Let $V$ be a module
over $k$ and let $S(V)$ be the symmetric tensor algebra over $V$;
$S(V)= k \oplus V
\oplus S^2(V) \oplus ...$. If $V$ is free of dimension $n$ generated by
$x_1,...,x_n$ then $S(V)$ is the polynomial algebra $k[x_1,...,x_n]$.
Assume that $V$ is a flat $k$-module (e.g. a free module). Then there
is an isomorphism
$$S(V)\otimes \Lambda^nV \cong  HH_n(S(V))$$
where $\Lambda^*V$ is the exterior algebra of $V$.
If $V$ is free of dimension $n$ generated by
$x_1,...,x_n$ then $\Lambda^mV$ is a free $n\choose m$ $k$-module with
a basis $v_{i_1} \wedge v_{i_2} \wedge ...\wedge v_{i_m}$, $i_1 < i_2 <...
<i_m$.
\end{theorem}
The above theorem is a special case of Hochschild-Konstant-Rosenberg 
theorem about Hochschild homology of smooth algebras \cite{HKR}, 
which we discuss in section 4. Here we stress, after Loday, that 
the isomorphism $\varepsilon_*: S(V)\otimes \Lambda^*V \to  HH_*(S(V))$
is induced by a chain map, that is not true in general for smooth algebras.
$S(V)\otimes \Lambda^*V$ is a chain complex with the zero boundary maps.
The chain map $\varepsilon: S(V)\otimes \Lambda^*V \to C_n(S(V),S(V))$ 
is given by $\varepsilon (a_0\otimes a_1 \wedge ...\wedge a_n) = 
\varepsilon_n (a_0,a_1,...,a_n)$ where $\varepsilon_n$ is the 
antisymmetrization map given as the sum 
$\sum_{\sigma \in S_n}sgn(\sigma)\sigma((a_0,a_1,...,a_n))$ and 
the permutation $\sigma \in S_n$ acts by $\sigma((a_0,a_1,...,a_n)) = 
(a_0,a_{\sigma^{-1}(1)},..., a_{\sigma^{-1}(n)})$.

In the second section we describe Khovanov homology of links and its 
(``comultiplication-free") version for graphs 
introduced in \cite{H-R-2}. In order 
to compare them with Hochschild homology we offer various generalizations,
relaxing a condition that underlying algebra needs to be commutative 
(in the case of a polygon or a line graph), and 
allowing an $\M$-reduced case. 
Another innovation in our presentation is that we start with 
cohomology of a functor from the category 
of subsets of a fixed set to the category of modules. In this setting we 
describe graph homology introduced in \cite{H-R-2} and its generalizations 
to $\mathbb M$-reduced cohomology and to cohomology of a polygon or 
a line graphs with a noncommutative underlying algebra. 
A set in this case is the set of edges of a graph. 
Finally, we describe a generalization 
to ``supersets" which allows homology for signed graphs, link diagrams 
and, in some cases, links (with the classical Khovanov 
homology as the main example). In examples of functors from ``supersets" 
we utilize comultiplication in the underlying algebra $\A$, after 
Khovanov \cite{Kh-1}, we assume that $\A$ is a Frobenius algebra 
\cite{Abr,Kock}.

In the third section we prove our main result relating 
 Hochschild homology $HH_*(\A)$ to graph cohomology and 
Khovanov homology of links, and the homology $H_n(\A,\M)$ 
to reduced ($\M$-reduced, 
more precisely)  cohomology of graphs and links.

In the fourth section we use our main result to describe graph cohomology 
of polygons for various algebras solving, in particular, several 
problems from \cite{H-P-R}.

We envision the connection between Connes cyclic homology and Khovanov 
type homology, possibly when analyzing symmetry of graphs and links.
For the idea of cyclic homology it is the best to quote from J-L. Loday
\cite{Lo}:\\ 
``...in his search for a non-commutative analogue of de Rham homology
theory, A. Connes discovered in 1981 the following striking phenomenon: \\
- the Hochschild boundary map $b$ is still well defined when one factors
out the module $A\otimes A^{\otimes n} = A^{\otimes n+1}$ by the action
of the (signed) cyclic permutation of order $n+1$.\\
Hence a new complex was born whose homology is now called (at least in
characteristic zero) {\it cyclic homology}."

\section{Khovanov homology}

We start this section with a very abstract definition based 
on Khovanov construction but, initially, devoid of topological 
or geometric context.
In this setting 
we recall and generalize the concept of graph cohomology \cite{H-R-2}
and of classical Khovanov 
homology for unoriented framed links. We follow, in part, the exposition in 
\cite{H-R-1,H-R-2} and \cite{Vi}.
We review, after \cite{H-P-R}, the connection between Khovanov homology 
of links and graph cohomology of associated Tait graphs.
We define homology of link diagrams 
related to graph cohomology for any commutative algebra\footnote{We can extend 
the definition to noncommutative algebras in the case of $(2,n)$-torus 
link diagrams.}.

\subsection{Cohomology of a functor on sets}\label{sub2.1}\ \\ 
Let $k$ be a commutative ring and $E$ a finite (or countable) set.
\begin{definition}
Let $\Phi$ be a functor from the category of subsets of $E$
(i.e. subsets of $E$ are
objects and inclusions are morphisms) to the category of $k$-modules.
We define the ``Khovanov cohomology" of the functor, $H^i(\Phi)$,
 as follows. 
We start from the graded $k$-module $\{C^i(\Phi)\}$, where $C^i(\Phi)$ 
is the direct sum of $\Phi(s)$ over all
$s\in E$ of $i$ elements ($|s|=i$). 
To define $d: C^i(\Phi) \to C^{i+1}(\Phi)$ we first define 
face maps $d_e(s)= \Phi(s\subset s\cup e)(s)$ for $e \in E-s$ (notice that 
$s\subset s\cup e$ is the unique morphism in $Mor(s,s\cup e)$). 
Now, as usually, $d(s)= \sum_{e \notin s} (-1)^{t(s,e)}d_e(s)$,
where $t(s,e)$ requires ordering of elements of $E$ and is equal the
number of elements of $s$ smaller then $e$. Because $\Phi$ is a functor, 
therefore we have $d_{e_2}d_{e_1}(s) = d_{e_1}d_{e_2}(s)$ for any 
$e_1,e_2 \notin s$. The sign convention guarantees now that $d^2=0$ and 
$(\{C^i(\Phi)\},d)$ is a cochain complex.
Now we define, in a standard way, cohomology as
 $H^i(\Phi)= ker (d(C^i(\Phi)\to C^{i+1}(\Phi))/d(C^{i-1}(\Phi))$.
The standard argument shows that $H^i(\Phi)$ is independent on 
ordering of $E$.
\end{definition}
In the case when $E$ are edges of a graph $G$ we can define
specific functors in various ways taking into account a structure of $G$. 
We construct below our main example: a generalization 
of a graph cohomology, defined in \cite{H-R-2}, to $\M$-reduced case 
and its translation to 
homology of alternating diagrams. In the case of the 
algebra $\A_2=Z[x]/(x^2)$,
this homology agrees partially with the classical Khovanov homology (see 
Theorem 2.7). To deal with all link diagrams we will later  
expand Definition 2.1 to ``supersets" (Definition 2.4) in a construction 
which can incorporate multiplication and comultiplication 
in $\A$ (Example 2.5).

\begin{definition}[Cohomology introduced in \cite{H-R-2} and extended to 
$\M$-reduced case]\ \
\begin{enumerate}
\item[(1)] We define here $\M$-reduced cohomology denoted by by 
$H^*_{\A,\mathbb M}(G,v_1)$.
If we assume $\mathbb M = \A$ we obtain
(comultiplication-free) cohomology of graphs, $H^*_{\A}(G)$,
  defined in \cite{H-R-2}.\\
Let $G$ be a graph with an edge set $E=E(G)$, and a chosen base vertex $v_1$.
Fix a commutative $k$-algebra $\A$ and an $\A$-module $\mathbb M$. 
We define a functor $\Phi$ on a category of subsets of $E$ as follows:
\begin{enumerate}
\item[(Objects)] To define the functor $\Phi$ on objects $s\subset E$ we 
define it more generally on any subgraph $H \subset G$, starting 
from a connected $H$, to be
$\Phi(H)=\mathbb M$ if $v_1 \in H$, and $\Phi(H)=\A$ if $v_1 \notin H$.
If $H$ has connected components $H_1,...,H_k$ then we define 
$\Phi(H)= \Phi(H_1)\otimes ... \otimes \Phi(H_k)$. Finally, 
$\Phi(s)= \Phi([G:s])$ where $[G:s]$ is a subgraph of $G$ containing 
all vertices of $G$ and edges $s$. \\
In what follows $k(s)$ is the number of components of $[G:s]$.
\item[(Morphisms)] It suffices to define $\Phi(s \subset s\cup e)$ 
where $e\notin s$. The definition depends now on the role of $e$ in 
$[G:s]$ as follows:
\begin{enumerate}
\item[(i)] Assume that $e$ connects different components of $[G:s]$:\\
(i') If $e$ connects components $u_i$ and $u_{i+1}$ not containing $v_1$ 
then we define 
$\Phi(s \subset s\cup e)(m,a_1,...,a_i,a_{i+1},...,a_{k(s)-1})=
(m,a_1,...,a_ia_{i+1},...,a_{k(s)-1})$.\\
(i'') If $e$ connects a component of $[G:s]$ containing $v_1$ with another 
component of $[G:s]$, say $u_1$, then we put 
$\Phi(s \subset s\cup e)(m,a_1,...,a_{k(s)-1})= (ma_1,...,a_{k(s)-1})$.
\item[(ii)] Assume that $e$ connects vertices of the same component 
of $[G:s]$, then $\Phi(s \subset s\cup e)$ is the identity map 
on $\Phi(s)=\Phi(s\cup e)= \mathbb M \otimes \A^{\otimes k(s)-1}$ 
\end{enumerate}
In the proof that $\Phi$ is a functor commutativity of $\A$ is important 
(compare (3)). 
\end{enumerate}
\item[(2)] If we modify the functor $\Phi$ from (1) to a new functor,
 $\hat\Phi$ which differs from $\Phi$ only in the rule $(1)(ii)$ 
that $\hat\Phi(s \subset s\cup e)$ is a zero map if $e$ 
connects vertices of the same component of $[G:s]$. 
We denote the cohomology yielded by the functor $\hat\Phi$ by 
$\hat H^*_{\A,\mathbb M}(G,v_1)$.  
For $\mathbb M = \A$ this cohomology, $\hat H^*_{\A}(G)$,
 is introduced in \cite{H-R-2}.  
\item[(3)]  For a polygon or a line graph the cohomology 
$\hat H^*_{\A,\mathbb M}(G,v_1)$ is defined also for 
a noncommutative algebra $\A$ and any $\A$-bimodule $\M$. 
We consider a polygon or a line graph as a directed graph: from left 
to right (in the case of a line graph, Fig. 3.1) and in the anti-clockwise 
orientation (in the case of a polygon).
In the formula for the morphism, $\hat\Phi(s\subset s \cup e)$, of 
Definition 2.2(2)
we use the product $xy$ if $x$ is the weight of the initial point of the
directed edge $e$ connecting different components of $[G:s]$.\
 We use this graph cohomology of the directed polygon when
comparing graph cohomology with Hochschild homology.
\end{enumerate} 
\end{definition} 
Notice that cohomology described in (1) and (2) coincide to certain degree.
Namely $H^i_{\A,\mathbb M}(G,v_1) = \hat H^i_{\A,\mathbb M}(G,v_1)$ 
for all $i < \ell -1$ where $\ell$ is 
the length of the shortest cycle in $G$.\\
Furthermore, if $k$ is a principal ideal domain (e.g. $k=Z$)
 and $\A$ and $\M$ are free 
$k$-modules then $Tor(H^i_{\A,\mathbb M}(G,v_1)) = 
Tor(\hat H^i_{\A,\mathbb M}(G,v_1))$ for $i=\ell-1$ (compare Theorem 2.7).

\begin{remark} 
\begin{enumerate}
\item[(i)] One can generalize\footnote{We are motivated here by 
 \cite{Sto}.} construction in Definition 2.2(1) and (2)
 by choosing the sequence of elements $f_1,f_2,f_3,...,f_{|E|}$ 
in $\A$ and modifying functors $\Phi$ and $\hat \Phi$ on morphisms 
to get the functors $\Phi'$ and $\hat \Phi'$. 
We put $\Phi'(s \subset s\cup e)=f_{|s|+1}\Phi(s \subset s\cup e)$ 
and, similarly,  
$\hat\Phi'(s \subset s\cup e)=f_{|s|+1}\hat\Phi(s \subset s\cup e)$.  
\item[(ii)] If $\A$ is not commutative and we work with cohomology 
of a line graph or a polygon (as in Definition 2.2(3))
we have to assume, in order to have $d^2=0$,
 that $f_i$'s are in the center of $\A$. 
We define 
$f\cdot(m,a_1,...,a_{k(s)-1})(s)=(fm,a_1,...,a_{k(s)-1})(s)$.
\end{enumerate}
\end{remark}

We can define Khovanov cohomology on an alternating link diagram, $D$ 
by considering associated plane graph, G(D) (Tait graph; compare Fig. 2.1) 
and its cohomology described in Definition 2.2 and Remark 2.3. 

In order to define Khovanov cohomology on any link diagram (and 
take both multiplication and comultiplication into account) we 
have to define cohomology on any signed planar graph. 
 We can start, as in Definition  2.1, from the 
very general setting (again cohomology of a functor) and to produce 
 specific examples of a cohomology of  signed planar graphs using 
 a coherent algebra and coalgebra structures (Frobenius algebra).
\begin{definition}\label{2.4}
Let $k$ be a commutative ring and $E=E_+ \cup E_-$ a finite set divided 
into two disjoint subsets (positive and negative sets). We consider the 
category of subsets of $E$ 
($E\supset s=s_+\cup s_-$ where $s_{\pm}= s\cap E_{\pm}$). the set $Mor(s,s')$ 
is either empty or has one element if $s_- \subset s'_-$ and 
$s_+ \supset s'_-$. Objects are graded by $\sigma(s)=|s_-|- |s_+|$. 
Let us call this category the {\it superset category}
 (as the set $E$
 is initially $Z_2$-graded). We define ``Khovanov cohomology" for every functor,
$\Phi$, from the superset category to the category of $k$-modules.
We define cohomology of $\Phi$ in the similar way as for a functor 
from the category of sets (which corresponds to the case $E=E_-$).
The cochain complex corresponding to $\Phi$ is defined to be
$\{C^i(\Phi)\}$ where $C^i(\Phi)$ is the direct sum of $\Phi(s)$ over all
$s\in E$ with $\sigma(s)=i$. To define $d: C^i(\Phi) \to C^{i+1}(\Phi)$ 
we first define face maps $d_e(s)$ where $e=e_- \notin s_-$ ($e_- \in E_-$) 
or $e=e_+ \in s_+$. In such a case 
$d_{e_-}(s)= \Phi(s\subset s\cup e_-)(s)$  and 
$d_{e_+}(s)= \Phi(s\supset s-e_+)$. We define  
$d(s)= \sum_{e \notin s} (-1)^{t(s,e)}d_e(s)$,
where $t(s,e)$ requires ordering of elements of $E$ and is equal the
number of elements of $s_-$ smaller then $e$ plus the number of elements 
of $s_+$ bigger than $e$.
We obtain the cochain complex whose cohomology does not 
depend on ordering of $E$.
\end{definition}
\begin{example}\label{2.5}
Let $G$ be a signed plane graph with an edge set 
$E=E_+ \cup E_-$ were $E_+$ is the 
set of positive edges and $E_-$ is the set of negative edges. We  
define the functor from the superset category $\mathbb E$ using the 
fact that $G$ is the (signed) Tait graph of a link diagram $D(G)$ 
(with white infinite region). See Figure 2.1 for conventions).
\\ \ \\
\centerline{\psfig{figure=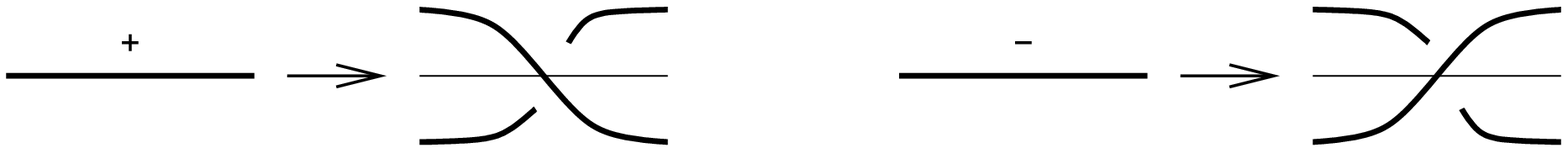,height=1.2cm}}
\centerline{ Figure 2.1} \ \\

To define the functor $\Phi$ we fix 
a Frobenius algebra $\A$ with multiplication $\mu$ and comultiplication 
$\Delta$ (the main example used by Khovanov
 is the algebra of truncated polynomials $\A_m=Z[x]/(x^m)$ with a 
coproduct $\Delta (x^k)=\sum_{i+j=m-1+k}x^i \otimes x^j$).
To get  our functor on objects, $\Phi(s)$,
 we consider the Kauffman state defined by $s$ (so also 
denoted by $s$) and we associate $\A$ to every circle of $D_s$ obtained 
from $D(G)$ by smoothing every crossing according to $s$ and 
then taking tensor product of these copies of $\A$ (compare \cite{Vi}).
To get $\Phi$ on morphisms, 
without going into details, we say succinctly that 
if Mor$(s,s')\neq \emptyset$ 
 and $\sigma(s')=\sigma(s)+1$ then $\Phi(Mor(s,s'))$ is defined 
using product or coproduct depending on whether $D_s$ has more or less 
circles than $D_{s'}$. For $A_2$ we obtain the classical Khovanov homology.  
\end{example}
\begin{remark} \label{2.6} One can also extend Example 2.5 to 
include the concept of $M$-reduced cohomology. We can consider, for example 
$M$ to be an ideal in $\A$ with $\Delta(\M) \subset M\otimes M$. An example, 
considered by Khovanov, is $A_m$ with $M$ generated by $x^{m-1}$ 
(in which case $\Delta(x^{m-1})= x^{m-1}\otimes x^{m-1}$). 
\end{remark}
One can build more delicate (co)homology theory for ribbon graphs 
(flat vertex graph) using the fact that they embed uniquely 
into the closed surface.
For $A=A_2$ it can be achieved using the approach presented 
in \cite{APS}, while for more general (Frobenius) algebras
it is not yet done (most likely one should not use Frobenius algebra 
alone but its proper enhancement like in $A_2$ case).

In \cite{H-P-R} we proved the following relation between 
graph cohomology and classical Khovanov homology 
of alternating links.
\begin{theorem}\label{Let $D$ be the diagram}
Let $D$ be the diagram of an unoriented framed alternating link
and let $G$ be its Tait graph. Let $\ell$ be the length of the
shortest cycle in $G$. For all $i<\ell-1$, we have

\[
H^{i,j}_{\A_2}(G) \cong H_{a,b}(D)
\]
with $\left \{\begin{array}{l}
a=E(G)-2i,\\
b=E(G)-2V(G)+4j.
\end{array}
\right.$

where $H_{a,b}(D)$ are the Khovanov homology groups 
of the unoriented framed link defined by $D$, as explained in
\cite{Vi}.

Furthermore, $Tor(H^{i,j}_{\A_2}(G)) = Tor(H_{a,b}(D))$ for $i=\ell-1$.
\end{theorem}

We also speculate that for other $sl(m)$ Khovanov-Rozansky 
homology \cite{Kh-R-1,Kh-R-1} the graph cohomology of 
a polygon with $\A=\A_m=Z[x]/(x^m)$ keeps 
essential information on $sl(m)$-homology of a torus link.

\section{Relation between Hochschild homology and Khovanov homology}

The main goal of our paper is to demonstrate relation between 
Khovanov homology and Hochschild homology. Initially I observed 
this connection 
for a commutative algebra $\A$ by showing that for every commutative
unital algebra $\A$ the graph cohomology of 
$(n+1)$-gon, $H^{i}_{\A}(P_{n+1})$, is isomorphic to Hochschild homology 
of $\A$, $H_{n-i}(\A)$; $0<i<n$. From this, via Theorem 2.7, relation 
between classical Khovanov homology of $(2,n+1)$ torus link and 
Hochschild homology of $\A_2=Z[x]/(x^2)$, follows. This relation was 
also observed independently by Magnus Jacobson \cite{Jac}.
  
In this paper  we prove more general result. In order to 
formulate it we use an extended version of (Khovanov type) 
graph cohomology (working with noncommutative algebras and 
${\mathbb M}$-reduced cohomology):  
\begin{enumerate}
\item[(i)] We can work with a noncommutative algebra, because, as mentioned 
in Section 2 (Definition 2.2(3)), for a polygon the graph cohomology 
is defined also for noncommutative algebras.

\item[(ii)] We can fix an $\A$-bialgebra ${\mathbb M}$ 
and compare ${\mathbb M}$-reduced Khovanov type graph cohomology with 
the Hochschild homology of ${\A}$ with coefficients in 
${\A}$-bimodule ${\mathbb M}$.
\end{enumerate}
One can observe that we generalize notion of graph (co)homology while 
we keep the original definition of Hochschild homology. Our point is 
that the graph (co)homology is the proper generalization of 
Hochschild homology: from a polygon to any graph. We have this 
interpretation only for a commutative $\A$. It seems to be, that if 
one work with general graphs and not necessary commutative algebras then 
these algebras should satisfy some ''multiface" properties. 
Very likely that 
planar algebras or operads provide the proper 
framework\footnote{The author's idea of working with directed graphs (quivers)
seems to work, as observed by Y.Rong, only for line graphs and polygons.}.

The interpretation of Hochschild homology as
a homology of $\A$ treated as an algebra over 
$\A\otimes \A^{op}$ allows us to use the standard tool of homological algebra,
that is  we 
find appropriate (partial) free resolution of $\A\otimes \A^{op}$ 
module $\A$ using graph cochain complex of a line graph (Figure 3.1). 
The graph cohomology of 
the polygon is the cohomology obtained from this resolution. In Theorem 3.1 
we use cohomology $ \hat H^i_{\A,\mathbb M}(P_{n+1})$ because we 
do not assume that $\A$ is commutative. 

\begin{theorem}\label{3.1} Let $\A$ be a unital 
algebra\footnote{We assume in this paper that $\A$ is a 
free $k$-module,
but we could relax the condition to have $\A$ to be projective or, 
more generally, flat over a commutative ring with identity $k$; 
compare \cite{Lo}. We require $\A$ to be a unital algebra in order to 
have an isomorphism $\M \otimes_{\A^{\epsilon}} \A^{\otimes n+2} = 
\M \otimes \A^{\otimes n}$; the isomorphism is given by 
$\M \otimes_{\A^{\epsilon}} \A^{\otimes n+2} \ni(m,a_0,a_1,...,a_n,a_{n+1}) 
\to (a_{n+1}ma_0,1,a_1,...,a_n,1)$ which we can write succinctly as 
$(a_{n+1}ma_0,a_1,...,a_n) \in \M \otimes \A^{\otimes n}$. We should stress 
that in $\M \otimes_{\A^{\epsilon}} \A^{\otimes n+2}$ the tensor product 
is taken over $\A^{\epsilon}=\A\otimes \A^{op}$ while in 
$\M \otimes \A^{\otimes n}$ the tensor product is taken over $k$.}, 
$\mathbb M$ an $\A$-bimodule and 
$P_{n+1}$ -- the $(n+1)-gon$. Then for $0< i \leq n$ we have:
$$ \hat H^i_{\A,\mathbb M}(P_{n+1}) = H_{n-i}(\A,\mathbb M).$$
Furthermore, if ${\A}$ is a graded algebra and ${\mathbb M}$ a 
coherently graded module then $ \hat H^{i,j}_{\A,\mathbb M}(P_{n+1}) = 
H_{n-i,j}(\A,\mathbb M)$, for $0< i \leq n$ and every $j$.
\end{theorem}
\begin{corollary}
$ \hat H^{i,j}_{\A}(P_{n+1}) =
HH_{n-i,j}(\A)$, for $0< i \leq n$ and every $j$.
Furthermore, for a commutative $\A$, 
$H^{i,j}_{\A}(P_{n+1}) = \hat H^{i,j}_{\A}(P_{n+1})$ , for $0< i <n$ and 
$H^{n,*}_{\A}(P_{n+1}) = 0$, $\hat H^{n,j}_{\A}(P_{n+1})=  
HH_{0,*}(\A)=\A$.\ 
For a general $\A$,  $\hat H^{n,*}_{\A}(P_{n+1}) =HH_{0,*}(\A)=\A/(ab-ba)$).
\end{corollary}
{\bf Poof of Theorem 3.1}\\
We consider graph cohomology for a unital, possibly noncommutative
algebra $\A$ and any $\A$-bimodule $\mathbb M$. 
There is no difference in the proof 
between commutative and noncommutative case except that we have to 
prove some property of cohomology given in \cite{H-R-2} for 
a commutative $\A$ (Lemma 3.3).\\
The main idea of our proof is to interpret the graph cochain complex 
of a line graph as a (partial) resolution of $\A$. It was proved 
in \cite{H-R-2} that $H^i_{\A}$(line graph)$=0$, $i>0$  
for a commutative algebra $\A$. We give here the proof 
for any unital algebra $\A$.
Let $L_n$ be the (directed) line graph of $n+1$ vertices ($v_0,...,v_n$)
 and $n$ edges ($e_1,...,e_n$); Fig. 3.1. 
\\ \ \\
\centerline{\psfig{figure=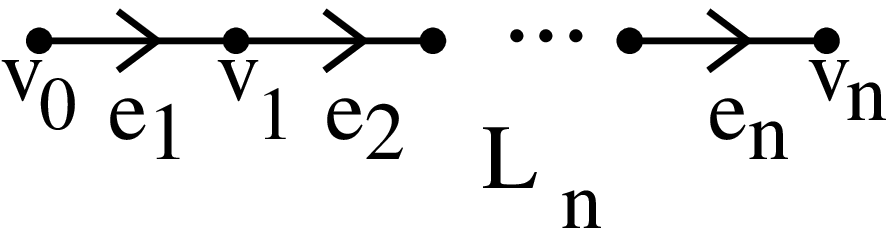,height=1.6cm}}
\centerline{ Figure 3.1} \ \\

\begin{lemma}
The graph cochain complex of $L_n$:\\
$C^*_{\A}(L_n):\ \   C^{0} \stackrel{d^0}{\to} C^{1} \stackrel{d^1}{\to} C^{2} 
...\to  C^{n-1} \stackrel{d^{(n-1)}}{\to} C^{n} $, 
is acyclic, except for the first term. 
That is, $\hat H^{i}_{\A}(L_{n})=0$ for $i>0$ 
and $\hat H^{0}_{\A}(L_{n})$ is usually nontrivial\footnote{From the 
fact that chromatic polynomial of $L_n$ is equal to $\lambda(\lambda -1)^n$ 
follows that  rank($\hat H^{0}_{\A}(L_{n})=$ 
rank($\A)$ $(rank(\A)-1)^n)$, we assume here that $k$ is a principal ideal 
domain. It was proven in \cite{H-R-2} that for 
a commutative $\A$ decomposable into 
$k1 \oplus \A/k$ one has that $H^*_{\A}(L_n)= H^0_{\A}(L_n)= 
\A \otimes (\A/k)^{\otimes n}$.}. 
\end{lemma}
For a line graph $\hat H_{\A} = H_{\A}$ so we will use $H_{\A}$ to simplify 
notation. \ 
We prove Lemma 3.3 by induction on $n$. For $n=0$, $L_0$ is 
the one vertex graph, thus $H^*_{\A}= H^0_{\A}= \A$ and 
the Lemma 3.3 holds. Assume that the lemma holds for $L_k$ with $k<n$.
In order to perform inductive step we construct the long exact sequence
of cohomology of line graphs. For a commutative ${\A}$ it is a special 
case\footnote{One can construct an exact sequence of functor cohomology 
imitating deleting-contracting exact sequence. One have to define properly 
two functors on $E\cup e$, one ``covariant" and one ``contravariant" and 
the exact sequence will be based on a functor on subsets of $E$ 
and these two additional functors. We will discuss this idea
 in a sequel paper.}
 of the exact sequence of graph cohomology in \cite{H-R-2} which in turn 
resemble the skein exact sequence of Khovanov homology \cite{Vi}:
$$ 0 \to H^{0}_{\A}(L_{n}) \to 
H^{0}_{\A}(L_{n})\otimes {\A} \stackrel{\partial}{\to}
H^{0}_{\A}(L_{n-1}) \to H^{1}_{\A}(L_{n}) \to...$$
$$\to H^{i-1}_{\A}(L_{n-1}) \to 
H^{i}_{\A}(L_{n}) \to 
H^{i}_{\A}(L_{n})\otimes {\A}\to...$$
such that $\partial:H^{0}_{\A}(L_{n}) \otimes {\A} \to 
H^{0}_{\A}(L_{n-1})$ is an epimorphism. 
From this exact sequence the inductive step follows.\\
To construct the above exact sequence we consider the short 
exact sequence of chain complexes (for the notation see Figure 3.1):
$$ 0 \to C^{i-1}(L_n/e_{n}) \stackrel{\alpha}{\to}  C^{i}(L_n) 
\stackrel{\beta}{\to} 
C^{i}(L_{n}-e_n) \to 0.$$
This exact sequence is constructed in the same way as in the case of 
commutative ${\A}$, that is ${\alpha}(a_0\otimes a_1\otimes ...
\otimes a_{n-i})(e_{j_1},...,e_{j_{i-1}})= (a_0\otimes a_1\otimes ...
\otimes a_{n-i})(e_{j_1},...,e_{j_{i-1}},e_n)$. Further, 
$\beta$ is defined in such a way that if $e_n\in s$ then $\beta (S)=0$,
and if  $e$ is not in $s$ then $\beta$ is the identity 
map (up to $(-1)^{|s|}$).


Exactness of the sequence 
follows from the definition.  This exact sequence leads to the long 
exact sequence of cohomology:
$$0 \to H^{0}_{\A}(L_{n}) \to
H^{0}_{\A}(L_{n}-e_n) \stackrel{\partial}{\to}
H^{0}_{\A}(L_{n}/e_n) \to $$
$$...\to H^{i-1}_{\A}(L_{n}/e_n) \to
H^{i}_{\A}(L_{n}) \to
H^{i}_{\A}(L_{n}-e_n) \to...$$
Now $L_{n-1}=L_{n}/e_n$ and $L_{n}-e_n$ is $L_{n-1}$ with an additional 
isolated vertex, therefore by a K\"unneth formula (see for example 
1.0.16 of \cite{Lo}) we have $H^{i}_{\A}(L_{n}-e_n)= 
H^{i}_{\A}(L_{n})\otimes {\A}$. From this we get 
the exact sequence used in the proof of Lemma 3.3. To see the 
epimorphism of $\partial$ notice that the map 
$H^{0}_{\A}(L_{n}-e_n) \to
H^{0}_{\A}(L_{n-1})$ is an epimorphism almost by the definition 
(we can think of decorating the last vertex of $L_n$ by $1$, to see 
the epimorphism. \\
There is chain map epimorphism $\partial_c: C^i_{\A}(L_{n}-e_n) \to
C^{i}_{\A}(L_{n-1})$ which is obtained by multiplying the 
weight of component containing $v_{n-1}$ by the weight of $v_n$ which 
descends to $\partial$. This map has a chain map section
$\partial_c^{-1}: C^{i}_{\A}(L_{n-1}) \to C^i_{\A}(L_{n}-e_n)$ 
such that in the image $v_n$ has always weight $1$. Because 
$\partial_c\partial_c^{-1} = Id$ therefore, on the 
cohomology level $\partial$ is an epimorphism.

We can continue now with the proof of Theorem 3.1. 
The (partially) acyclic chain complex of Lemma 3.3 is the chain complex 
of ${\A}^e={\A}\otimes {\A}^{op}$ modules. It is a 
(partial) free resolution of the ${\A}^e$-module ${\A}$. 
Upon tensoring this resolution with ${\mathbb M}$ considered as a right 
module over ${\A}^e$ we obtain the cochain complex, 
$$\{\mathbb M \otimes_{{\A}^e} C^i\}_{i=0}^{n-1} : \ 
\M \otimes_{{\A}^e} C^0  \stackrel{\partial^0}{\to}
\mathbb M \otimes_{{\A}^e} C^1 \stackrel{\partial^1}{\to} ...  \to
\mathbb M \otimes_{{\A}^e} C^{n-2} \stackrel{\partial^{n-1}}{\to}
\mathbb M \otimes_{{\A}^e} C^{n-1} \to 0$$
 whose cohomology 
(except possibly $H^0$) are the Hochschild homology of ${\A}$ 
with coefficients in ${\mathbb M}$ (compare for example \cite{Wei}).
Having in mind relation between  indexing we get that 
$H^i = H_{n-i}({\A},{\mathbb M})$ for $i>0$. 
To get exactly the chain complex of the ${\mathbb M}$-reduced
 (directed) graph homology of $P_{n}$, $\hat H^*_{{\A},{\mathbb M}}(P_n)$,
 we extend this chain complex to 
$$\{\mathbb M \otimes_{{\A}^e} C^i\}_{i=0}^{n}\ :\ \
\mathbb M \otimes_{{\A}^e} C^0  \stackrel{\partial^0}{\to}
\mathbb M \otimes_{{\A}^e} C^1 \stackrel{\partial^1}{\to} ...  \to
\mathbb M \otimes_{{\A}^e} C^{n-1} \stackrel{\partial^{n-1}}{\to}
\mathbb M \otimes_{{\A}^e} C^{n}=\mathbb M \to 0$$
where the homomorphism
$\partial^{n-1}$ is the zero map. 


To complete the proof of Theorem 3.1 we show that this complex is 
exactly the same as the cochain complex of the ${\mathbb M}$-reduced
 (directed) graph cohomology of $P_{n}$. We consider carefully the 
map $\mathbb M \otimes_{{\A}^e} C^j \stackrel{\partial^j}{\to} 
\mathbb M \otimes_{{\A}^e} C^{j+1}$. In the calculation 
we follow  the proof of Proposition 1.1.13 of \cite{Lo}.
The idea is to ``bend" the line graph $L_n$ to the polygon $P_n$ and 
show that it corresponds to tensoring, over ${\A}^e$ with 
$\mathbb M$; compare Figure 3.2.
\\ \ \\
\centerline{\psfig{figure=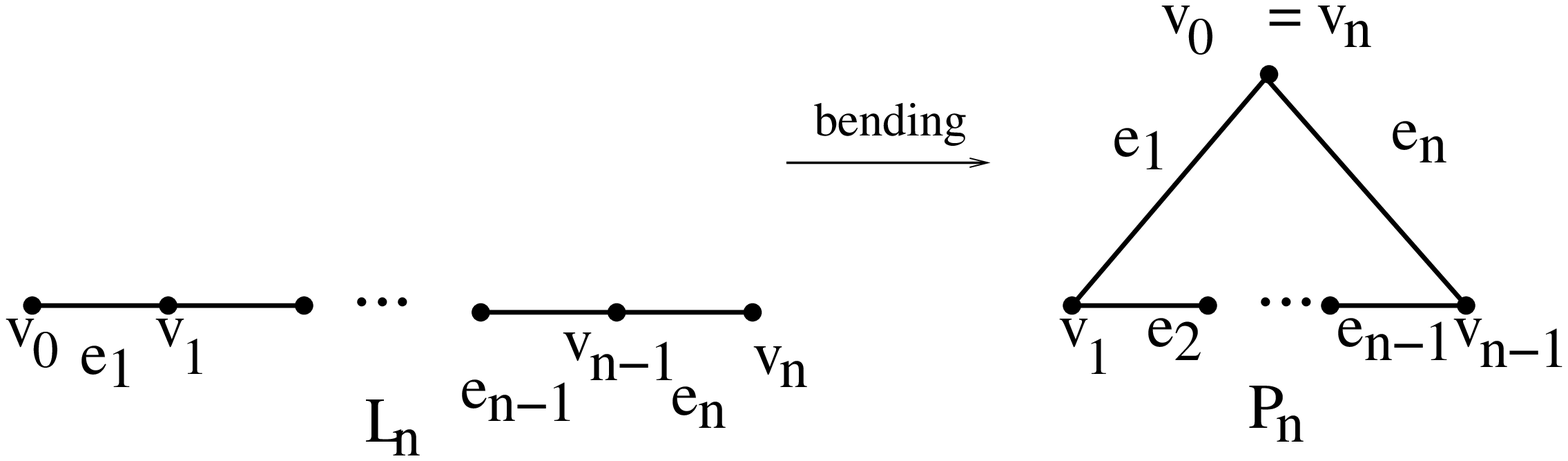,height=3.4cm}}
\centerline{ Figure 3.2} \ \\

Let us order components of $[G:s]$ ($G$ is equal to $P_n$ or $L_n$) 
in the anticlockwise orientation, 
starting from the component containing $v_0$ (decorated by an element of $\M$ 
if $G=P_n$). The element of $C^j_{\A}(L_n)$ will be denoted by 
$(a_{i_0},a_{i_1},...,a_{i_{n-j-1}},a_{i_{n-j}})(s)$ and of 
$C^j_{\A,\M}(P_n)$ by $(m,a_{i_1},...,a_{i_{n-j-1}})(s)$. In the 
isomorphism $\M \otimes_{{\A}^e} C^{j}_{\A}(L_n) \to C^j_{\A,\M}(P_n)$ 
($j<n$) the element $(m,a_{i_0},a_{i_1},...,a_{i_{n-j-1}},a_{i_{n-j}})(s)$ 
is send 
to $(a_{i_{n-j}}ma_{i_0},a_{i_1},...,a_{i_{n-j-1}})(s)$. 
One can easily check that it yields a cochain map (compare \cite{Lo}) so 
it induces the isomorphism on cohomology. Note that  
$(m,a_{i_0},a_{i_1},...,a_{i_{n-j-1}},a_{i_{n-j}}) = 
(a_{i_{n-j}}ma_{i_0},1,a_{i_1},...,a_{i_{n-j-1}},1)$ in 
$\M \otimes_{{\A}^e} C^{j}_{\A}(L_n)$.
The proof of Theorem 3.1 is completed.
 
\section{Calculations and speculations}

There is an extensive literature on Hochschild homology and a lot of 
ingenious methods of computing them (e.g. \cite{Lo,Wei,Ros,Kon}).
Our main result, Theorem 3.1, allows us to use these methods 
to compute graph cohomology  
for polygons and, to some extend, to other graphs (using for example an 
observation that some properties of a cohomology of a polygon propagate 
to graphs containing it (compare \cite{A-P,H-P-R}). Properties of 
Hochschild homology (and, equally well, cyclic homology) should 
eventually shed light on Khovanov type homology of links.

We start from adapting Theorem 1.3 about Hochschild homology of 
symmetric tensor algebra. The simplest case of one variable
polynomials $\A=\A_{\infty}=Z[x]$ allows us to extend Theorem 27 
of \cite{H-P-R} from the triangle to any polygon. 
$\A_{\infty}$ is a graded algebra with 
$x^i$ being of degree $i$. Consequently Hochschild homology of 
$\A_{\infty}$  is a bigraded
module. We treat $\A_{\infty}$ as a $Z$-module (an abelian group) and
to simplify description of homology we use the Poincar\'e polynomial
of $HH_{**}(\A_{\infty})$ to describe the free part of homology.
Recall that the Poincar\'e polynomial (or series) of
bigraded finitely generated $Z$-modules $H_{**}$
is $PP(t,q)=PP(H_{**})(t,p) =
\sum_{i,j}a_{i,j}t^iq^j$ where $a_{i,j}$ is the rank of the group $H_{i,j}$.
\begin{corollary}
For an $n$-gon $P_n$ the graph cohomology groups 
$H^{i,j}_{\A_{\infty}}(P_n)$ 
are free abelian with Poincar\'{e} polynomial $(q+q^2+q^3+...)^3 + 
t^{n-2}(q+q^2+q^3+...)=$ $ (\frac{q}{1-q})^3 + t^{n-2}\frac{q}{1-q}$.
\end{corollary}
\begin{proof} From Theorem 3.1 we obtain $H^{i,j}(P_n)$ for $0<i<n-1$.
 It was observed in \cite{H-R-2}  that $H^{i,j}(P_n)=0$ for $i \geq n-1$.
To find $H^{0,*}(P_n)$ we use the fact that the chromatic polynomial 
of $P_n$ is equal to $(\lambda -1)^n + (-1)^n(\lambda -1)$ and 
the graph cohomology categorify the chromatic polynomial. That is, 
if we substitute $t=-1$ and $1+q+q^2+q^3+...=\lambda$ in the Poincar\'{e} 
polynomial we obtain the chromatic polynomial \cite{H-R-2}.
\end{proof}
Another illustration of the power of our connection is for the 
algebra ${\A}={\A}_{p(x)}= Z[x]/(p(x))$ where $p(x)$ is a polynomial in 
$Z[x]$. We discuss the general case later, here let us notice that two 
special cases of $p(x)=x^m$ and $p(x)=x^m-1$ are of great interest in 
knot theory (in Khovanov-Rozansky homology \cite{Kh-R-1} 
and its deformations \cite{Gor}).
Let us apply first the knowledge of Hochschild homology for 
$\A_m = Z[x]/(x^m)$ (compare \cite{Lo}) and solving Conjectures 30 and 31 
of \cite{H-P-R}.
\begin{theorem}\label{4.2}\ \\ 
(Free)\ \ The Poincar\'e polynomial of $HH_{**}(\A_m)$ is equal to:\\
$(1+q+...+q^{m-1}) + t(q+q^2+...+q^{m-1}) + 
(t^2+t^3)(q+q^2+...+q^{m-1})q^{m} + (t^4+t^5)(q+q^2+...+q^{m-1})q^{2m}+...+ 
(t^{2i}+t^{2i+1})(q+q^2+...+q^{m-1})q^{im}+...$.
\ \  (Torsion)\ \
$Tor (H_{**}(\A_m))= \bigoplus_{i=1}^{\infty}H_{2i-1,im}(\A_m)$ 
where each summand is isomorphic to $Z_m$.
\end{theorem}
 We solve Conjectures 30 and 31 
 of \cite{H-P-R} by applying Theorems 4.2 and 3.1.
\begin{corollary}\label{4.3}
\begin{enumerate}
\item[(Odd)] For $n=2g+1$ we have: \\
$Tor(H^{*,*}_{{\A}_{m}}(P_{2g+1}))=
 H^{v-2,m}_{{\A}_{m}}(P_{2g+1})\oplus
H^{v-4,2m}_{{\A}_{m}}(P_{2g+1})\oplus ... \oplus
H^{1,gm}_{{\A}_{m}}(P_{2g+1})$ with each summand
isomorphic to $Z_m$.\\
The Poincar\'{e} polynomial of $H^{*,*}_{{\A}_{m}}(P_{2q+1})$ is
equal to\\
$ (q+...+q^{m-1})^v +$ \\
$(q+...+q^{m-1})(t^{v-2}+ (t^{v-3}+t^{v-4}q^m +
(t^2+t)q^{m(g-1})$. 
\item[(Even)]
For $n=2g+2$ we have: \\
$Tor(H^{*,*}_{{\A}_{m}}(P_{2g+2}))=
H^{v-2,m}_{{\A}_{m}}(P_{2g+1})\oplus
H^{v-4,2m}_{{\A}_{m}}(P_{2g+1})\oplus ... \oplus
H^{2,gm}_{{\A}_{m}}(P_{2g+2})$ with each summand
isomorphic to $Z_m$.\\
The Poincar\'{e} polynomial of $H^{*,*}_{{\A}_{m}}(P_{2q+2})$ is
equal to\\
$(q+...+q^{m-1})^n + q^{m(n/2)-1}(q+...q^{m-1})+$ \\
$(q+...+q^{m-1})(t^{n-2}+ (t^{n-3}+t^{n-4})q^m +
(t^3+t^2)q^{m(g-1)}+ tq^{mg})$.
\end{enumerate}
\end{corollary}
Assume that $m=2$ in Corollary 4.3. Then, using Theorem 2.7 we can 
recover Khovanov computation of homology of the torus 
link $T_{2,n}$ \cite{Kh-1,Kh-2}. In particular we get:
\begin{corollary}\cite{Kh-2}\  
Let $T_{2,-n}$ be a left-handed 
torus link  of type $(2,-n)$, $n>2$. Then the torsion part of the 
Khovanov homology of $T_{2,-n}$ is given by (in the description 
of homology we use notation of \cite{Vi} treating $T_{2,-n}$ as 
a framed link):\\  
(Odd)\ \  For $n$ odd, all the torsion of $H_{**}(T_{2,-n})$ is 
supported by \\ \ \  
$H_{n-2,3n-4}(T_{2,-n})= H_{n-4,3n-8}(T_{2,-n})=...= 
H_{-n+4,-n+8}(T_{2,-n})=Z_2$.\\
(Even)\ \ For $n$ even, all the torsion of $H_{**}(T_{2,-n})$ is
supported by\\
\ \ $H_{n-4,3n-8}(T_{2,-n})= H_{n-6,3n-12}(T_{2,-n})=...=
H_{-n+4,-n+8}(T_{2,-n})=Z_2$.\\

For a right-handed torus link  of type $(2,n)$, $n>2$, we can  
use the formula for the mirror image (Khovanov duality theorem; 
see for example\\
 \cite{A-P,APS}:\ \ \  $H_{-i,-j}({\bar D})= 
H_{ij}(D)/Tor(H_{ij}(D)) \oplus Tor(H_{i-2,j}(D))$.
\end{corollary}

The result on Hochschild homology of symmetric algebras has a 
major generalization to the large class of algebras called 
{\it smooth algebras}.
\begin{theorem}(\cite{Lo,HKR})\label{4.1}
For any smooth algebra $\A$ over $k$, the antisymmetrization map 
$\varepsilon_*: \Lambda^*_{\A|k} \to HH_n(\A)$ is an isomorphism of 
graded algebras. Here $\Omega^n_{\A|k}= \Lambda^n \Omega^1_{\A|k}$ 
is an $\A$-module of differential $n$-forms.
\end{theorem}
We refer to \cite{Lo} for a precise definition of a smooth algebra, 
here we only 
recall that the following are examples of smooth
algebras:\\
 (i) Any finite extension of a perfect field $k$ (e.g. a field
of characteristic zero). \\
(ii) The ring of algebraic functions on a nonsingular
variety over an algebraically closed field $k$, e.g. $k[x]$, $k[x_1,...,x_n]$,
$k[x,y,z,t]/(xt-yz-1)$ \cite{Lo}. \\ 
Not every quotient of a polynomial algebra is a smooth algebra.
For example, $C[x,y])(x^2y^3)$ or $Z[x]/(x^m)$ are not smooth.
The broadest, to my knowledge, treatment of Hochschild homology of 
algebras $C[x_1,...,x_n]/(Ideal)$ is given by Kontsevich in \cite{Kon}.
For us the motivation came from one variable polynomials, Theorem 40 
of \cite{H-P-R}. In particular we generalize Theorem 40(i) from a 
triangle to any polygon that is we compute the graph cohomology 
of a polygon for truncated
polynomial algebras and their deformations. Thus, possibly, we can 
approximate Khovanov-Rozansky $sl(n)$ homology and their deformations. 

\begin{theorem}\label{4.6}
\begin{enumerate}
\item[(i)] $HH_i({\A}_{p(x)}) = Z[x]/(p(x),p'(x))$ for $i$ odd and 
$HH_i({\A}_{p(x)}) = \{[q(x)]\in Z[x]/(p(x)) \ | \ \ 
q(x)p'(x)$  is divisible by $p(x)$ \}, 
for $i$ even $i \geq 2$. In both cases the $Z$ rank of the group 
is equal to the degree of gcd$(p(x),p'(x))$.
\item[(ii)] In particular, for $p(x)= x^{m+1}$, we obtain homology of 
the ring of truncated polynomials, ${\A}_{m+1}=Z[x]/x^{m+1}$ 
for which:\  \ 
$HH_i({\A}_{m+1})=Z_{m+1} \oplus Z^m$ for odd\ $i$ and 
$HH_i({\A}_{m+1})= Z^m$ for for $i$ even, $i\geq 2$. 
$HH_0({\A}_{m+1})= \A = Z^{m+1}$.
\item[(iii)] The graph cohomology of a polygon $P_n$, $H^i_{{\A}_{p(x)}}(P_n)$, 
is given by:\\
 $H^{n-2i}_{{\A}_{p(x)}}(P_n)= {\A}_{p(x)}/(p'(x))$ for $1\leq i \leq \frac{v-1}{2}$,\ 
and \\
$H^{n-2i-1}_{{\A}_{p(x)}}(P_n)= ker({\A}_{p(x)} 
\stackrel{p'(x)}{\to} {\A}_{p(x)})$ 
for $1\leq i \leq \frac{v-2}{2}$.\\
Furthermore,  $H^{k}_{{\A}_{p(x)}}(P_n)= 0$ for $k\geq n-1$ and 
$H^{0}_{{\A}_{p(x)}}(P_n)$ is a free abelian group of rank 
$(d-1)^n + (-1)^n(d-1)$ for $n$ even ($d$ denotes the degree of $p(x)$) and
  it is of rank $(d-1)^n + (-1)^n(d-1) - rank(H^{1}_{{\A}_{p(x)}}(P_n))$  
if $n$ is odd (notice that $(d-1)^n + (-1)^n(d-1)$ 
is the Euler characteristic of $\{H^{i}_{{\A}_{p(x)}}(P_n)\}$). 
 
\end{enumerate}
\end{theorem}

\begin{proof} Theorem 4.6(i)
 is proven by considering a resolution of ${\A}_{p(x)}$ as 
an ${\A}_{p(x)}^e={\A}_{p(x)}\otimes 
{\A}^{op}_{p(x)}$ 
module:
$$ \cdots \to {\A}_{p(x)}\otimes {\A}_{p(x)} 
\stackrel{u}{\to} 
{\A}_{p(x)}\otimes {\A}_{p(x)} \stackrel{v}{\to}
{\A}_{p(x)}\otimes {\A}_{p(x)} \stackrel{u}{\to} \cdots
\to {\A}_{p(x)} $$ 
where $u= x\otimes 1 - 1\otimes x$ and $v=\Delta (p(x))$ is a coproduct 
given by $\Delta (x^{i+1}) = x^i\otimes 1 + x^{i-1}\otimes x +...+ 
x \otimes x^{i-1} + 1\otimes x^i$.
\end{proof}
 Curious but not accidental observation 
is that by choosing coproduct $\Delta (1)= v$ 
we define a Frobenius algebra structure on $\A$. 
In Frobenius algebra $(x\otimes 1)\Delta (1) = (1 \otimes x)\Delta (1)$ 
which makes $uv=vu=0$ in our resolution. Furthermore the distinguished 
element of the Frobenius algebra $\mu\Delta (1)= p'(x)$. 
\\
\ \\
\ \\
\ \\


\ \\ \ \\ \ \\
\noindent \textsc{Dept. of Mathematics, Old Main Bldg., 1922 F St. NW
The George Washington University, Washington, DC 20052}\\
e-mail: {\tt przytyck@gwu.edu}

\begin{thebibliography}{99}

\bibitem[Abr]{Abr}
L. Abrams, Two-dimensional topological quantum field theories 
and Frobenius algebras, {\it Journal of the Knot Theory and Its 
Ramifications}, 5(5), 1995, 569-587.

\bibitem[A-P]{A-P}
M.~M. Asaeda, J.~H. Przytycki,
Khovanov homology: torsion and thickness,
{\it Advances in Topological Quantum Field Theory}, in
    Advances in topological quantum field theory, 135--166, Kluwer Acad.
    Publ., Dordrecht, 2004;
    e-print: {\tt http://www.arxiv.org/math.GT/0402402}.


\bibitem[APS]{APS}
M.~M. Asaeda,  J.~H. Przytycki, A.~S. Sikora, 
Categorification of the Kauffman bracket skein module of $I$-bundles
over surfaces,
{\it Algebraic \& Geometric Topology (AGT)}, 4, 2004, 1177-1210.
e-print: {\tt http://front.math.ucdavis.edu/math.QA/0403527}.

\bibitem [C-G-V]{C-G-V}
G. Cortinas, J. Guccione, O. Villamayor,
Cyclic homology of $K[Z/pZ]$, {\it K-theory}, 2, 1989, 603-616.

\bibitem [D-G-R]{D-G-R}
N. M. Dunfield, S. Gukov, J.Rasmussen, 
The Superpolynomial for Knot Homologies,
e-print:\ 
http://front.math.ucdavis.edu/math.GT/0505662

\bibitem [Gor]{Gor}
B. Gornik, Note on Khovanov link cohomology, e-print:\\
http://front.math.ucdavis.edu/math.QA/0402266

\bibitem [Gu-Sch-Va]{Gu-Sch-Va}
S. Gukov, A. Schwarz, C. Vafa,
Khovanov-Rozansky Homology and Topological Strings; \ e-print:\
http://arxiv.org/abs/hep-th/0412243

\bibitem [H-P-R]{H-P-R}
L. Helme-Guizon, J.~H. Przytycki, Y. Rong, Torsion in Graph Homology,\\
for Fundamenta Mathematicae, preprint (July 2005).\\
e-print: http://arxiv.org/abs/math.GT/0507245

\bibitem [H-R-1]{H-R-1}
L. Helme-Guizon, Y. Rong,
A Categorification for the Chromatic Polynomial,
to appear in {\it  Algebraic and Geometric Topology (AGT)}.\\
e-print:\ http://front.math.ucdavis.edu/math.CO/0412264

\bibitem [H-R-2]{H-R-2}
L. Helme-Guizon, Y. Rong, Graph Cohomologies from Arbitrary Algebras, 
\\ e-print:\ 
http://front.math.ucdavis.edu/math.QA/0506023

\bibitem [Hoch]{Hoch}
G. Hochschild, On the cohomology groups of an associative algebra, 
{\it Annals of Math.}, 46, 1945, 58-67.

\bibitem [HKR]{HKR}
G. Hochschild, B. Kostant, A. Rosenberg, Differential forms on regular
affine algebras, {\it Trans. Amer. Math. Soc.} 102, 1962, 383-408.

\bibitem [Jac]{Jac}
M. Jacobsson, Personal communication at 
AMS-IMS-SIAM Joint Summer Research Conference;
Quantum Topology--Contemporary Issues and Perspectives,
Snowbird, Utah,  June 5-9, 2005.

\bibitem[Kh-1] {Kh-1}
M.~Khovanov,
A categorification of the Jones polynomial,
{\it Duke Math. J.} 101 (2000), no. 3, 359--426,\ \ \
http://xxx.lanl.gov/abs/math.QA/9908171

\bibitem[Kh-2] {Kh-2}
M.~Khovanov, Patterns in knot cohomology I,
{\it Experiment. Math.} 12(3), 2003, 365-374, \ \ \
http://arxiv.org/abs/math/0201306

\bibitem[Kh-R-1]{Kh-R-1}
M. Khovanov, L. Rozansky, {\sl Matrix factorizations and link homology},
  e-print: {\tt http://www.arxiv.org/math.QA/0401268}.

\bibitem[Kh-R-2]{Kh-R-2}
M. Khovanov,  L. Rozansky, 
{\sl Matrix factorizations and link homology II},
   e-print: {\tt http://www.arxiv.org/math.QA/0505056}.

\bibitem [Kock]{Kock}
J. Kock, Frobenius algebras and 2D topological quantum field
theories, in {\it London Mathematical Society Student Texts}, 59,
Cambridge University Press, 2003.

\bibitem[Kon]{Kon} M. Kontsevich, Hochschild and Harrison 
cohomology of complete intersections, Appendix to the paper 
Quantization on Curves by C.Fronsdal, e-print:\\
http://front.math.ucdavis.edu/math-ph/0507021

\bibitem[Kon-2]{Kon-2} M. Kontsevich, 
Operads and Motives in Deformation Quantization, e-print:\\
http://front.math.ucdavis.edu/math.QA/9904055

\bibitem [Lo]{Lo} 
J-L. Loday, Cyclic Homology, Grund. Math. Wissen. 
Band 301, Springer-Verlag, Berlin, 1992 (second edition, 1998).

\bibitem [Pr]{Pr}
J. H. Przytycki, {\bf KNOTS:} From combinatorics of knot diagrams to the
combinatorial topology based on knots, Cambridge University Press,
accepted for publication, to appear 2007, pp. 600.

\bibitem [Ros]{Ros}
J. Rosenberg, Algebraic K-theory and its applications, 
Graduate Texts in Mathematics, 147, 1st ed 1994. Corr. 2nd printing, 1996,
X, 392 p.

\bibitem[Shu]{Shu}
A.~Shumakovitch, Torsion of the Khovanov Homology,
{\it Geometry and Topology (GT)}, to appear.
e-print:\  http://arxiv.org/abs/math.GT/0405474

\bibitem[Sto]{Sto}
M. Stosic, Categorification of the Dichromatic Polynomial for Graphs,
e-print:\  http://arxiv.org/abs/math.GT/0504239

\bibitem [Vi]{Vi}
O. Viro, Remarks on definition of Khovanov homology, e-print:\\
http://front.math.ucdavis.edu/math.GT/0202199

\bibitem [Wei]{Wei}
C. A. Weibel, An introduction to homological algebra, Cambridge studies 
in advanced mathematics, 38, Cambridge University Press,  1995.
\end{thebibliography}
\end{document}